\documentclass[a4paper,10pt]{article}
\usepackage[latin1]{inputenc}
\usepackage[T1]{fontenc}
\usepackage{amsmath}
\usepackage{amsfonts}
\usepackage{amstext}
\usepackage{amsthm}
\usepackage[all]{xy}
\usepackage{geometry}
\usepackage{srcltx}
\usepackage{graphics}
\usepackage{epsfig}
\usepackage{subfigure}
\usepackage{slashed}
\usepackage{color}
\usepackage{amssymb}
\usepackage{srcltx}
\newtheorem{theorem}{Theorem}[section]

\newtheorem{cor}[theorem]{Corollary}

\DeclareMathOperator{\sgn}{sgn}

\DeclareMathOperator{\diverg}{div}

\title{On bifurcation for semilinear elliptic Dirichlet problems and the Morse-Smale index theorem}
\author{Alessandro Portaluri and Nils Waterstraat}

\begin{document}
\date{}
\maketitle

\footnotetext[2]{A. Portaluri was supported by the grant PRIN2009 ``Critical Point Theory and Perturbative Methods
for Nonlinear Differential Equations''.}
\footnotetext[3]{N. Waterstraat was supported by a postdoctoral fellowship of
the German Academic Exchange Service (DAAD).}

\begin{abstract}
We study bifurcation from a branch of trivial solutions of semilinear elliptic Dirichlet boundary value problems on star-shaped domains, where the bifurcation parameter is introduced by shrinking the domain. In the proof of our main theorem we obtain in addition a special case of an index theorem due to S. Smale.
\end{abstract}

\section{Introduction}
Let $\Omega\subset\mathbb{R}^N$ be a smooth domain which is star-shaped with respect to $0$ and set

\begin{align*}
\Omega_r:=\{r\cdot x\in\mathbb{R}^N:\,\,x\in\Omega\},\quad 0<r\leq 1.
\end{align*}
We consider for $0<r\leq 1$ the semilinear elliptic Dirichlet problems

\begin{equation}\label{equation}
\left\{
\begin{aligned}
-\Delta u(x)+g(x,u(x))&=0,\quad x\in\Omega_r\\
u(x)&=0,\quad x\in\partial\Omega_r,
\end{aligned}
\right.
\end{equation}
where $g\in C^\infty(\overline{\Omega}\times\mathbb{R})$ and $g(x,0)=0$ for all $x\in\Omega$. Moreover, we assume that there exists $C>0$ such that $g$ satisfies the standard growth conditions

\begin{align*}
|g(x,\xi)|\leq C(1+|\xi|^\alpha),\quad \left|\frac{\partial g}{\partial\xi}(x,\xi)\right|\leq C(1+|\xi|^\beta),\quad (x,\xi)\in \Omega\times\mathbb{R}, 
\end{align*}
for certain constants $\alpha,\beta\geq 0$ depending on the dimension $N$ (cf. \cite[\S 1.2]{Ambrosetti}).\\
We call $r^\ast\in(0,1]$ a \textit{bifurcation point} for the equations \eqref{equation} if 
there exists a sequence $r_n\rightarrow r^\ast$ and $u_n\in H^1_0(\Omega_{r_n})$ such that $u_n$ is a non-trivial weak solution of the boundary value problem \eqref{equation} on $\Omega_{r_n}$ and $\|u_n\|_{H^1_0(\Omega_{r_n})}\rightarrow 0$.\\
Now denote $f(x)=\frac{\partial g}{\partial \xi}(x,0)$, $x\in\overline{\Omega}$, and consider the linearised equations

\begin{equation}\label{equationII}
\left\{
\begin{aligned}
-\Delta u(x)+f(x)u(x)&=0,\quad x\in\Omega_r\\
u(x)&=0,\quad x\in\partial\Omega_r.
\end{aligned}
\right.
\end{equation}
We call $r^\ast\in(0,1]$ a \textit{conjugate point} for \eqref{equationII} if

\begin{align*}
m(r^\ast):=\dim\{u\in C^2(\Omega_{r^\ast}):\, u\quad\text{solves}\quad\eqref{equationII}\}>0
\end{align*}
and henceforth we assume that $m(1)=0$. Our main theorem reads as follows:

\begin{theorem}\label{theorem}
The bifurcation points of \eqref{equation} are precisely the conjugate points of \eqref{equationII}. 
\end{theorem}

Our proof of Theorem \ref{theorem} uses crossing forms and the connection between bifurcation of branches of critical points and the spectral flow from \cite{SFLPejsachowicz}. 
In particular, we do not use any domain monotonicity properties of eigenvalues of \eqref{equationII}.\\
Interestingly, we moreover obtain from our computations a new and simple proof of a classical theorem due to Smale \cite{Smale} for the equations \eqref{equationII}. 
In what follows we denote by $M$ the Morse index of \eqref{equationII} on the full domain $\Omega_1=\Omega$, i.e. the number of negative eigenvalues counted according to their multiplicities.

\begin{cor}\label{corI}
$m(r)=0$ for almost all $0<r<1$ and 

\begin{align*}
M=\sum_{0<r<1}{m(r)}.
\end{align*} 
\end{cor}
We want to point out that Smale considered in \cite{Smale} the Dirichlet problem for general strongly elliptic differential operators on vector bundles over compact manifolds with boundary. 
It seems to us that a computation of the corresponding crossing forms, as we will do for the equation \eqref{equationII} in Section \ref{computations}, 
is no longer possible in this generality. However, Smale's theorem has attracted some interest in recent years for scalar equations on star-shaped domains in Euclidean spaces (cf. \cite{Chris} and \cite{Ale}).\\
Finally, from Theorem \ref{theorem} and Corollary \ref{corI} we immediately obtain the following result:

\begin{cor}\label{corII}
If $M\neq 0$, then there exist at least
	
	\begin{align*}
	\left\lfloor\frac{M}{\max_{0<r<1}m(r)}\right\rfloor
	\end{align*}
distinct bifurcation points in $(0,1)$.
\end{cor}
  
In the special case $N=1$, i.e. if \eqref{equation} is a semilinear ODE on a compact interval, we have $0\leq m(r)\leq 1$ and hence the number of bifurcation points in Corollary \ref{corII} is precisely the Morse index $M$.

%%%%%%%%%%%%%%%%%%%%%%%%%%%%%%%%%%%%%%%%%%%%%%%%%%%%%%%%%%%%%%%%%%%%%%%%%%%%%%%%%%%%%%%%%%%%%%%%%%%%%%%%%%%%%%%%%%%%%%%%%%%%%%%%%%%%%%%%%%%%%%%%%%%%%%%%%%%%%%%%%%%%%%%%%%%%%%%%%%%%%%%%%%%%%%%%%%%%%%%%%%%%%%%%%%%%%%%%%%%%%%%%%%%%%%%%%%%%%%%%%%%%%%%%%%%%%%%%%%%%%%%%%%%%%%%%%%%%%%%%%%%%%%%%%%%%%%%%%%%%%%%%%%%%%%%%%%%%%%%%%%%%%%%%%%%%%%%%%%%%%%%%%%%%%%%%%%%%%%%%%%%%%%%%%%%%%%%%%%%%%%%%%%%%%%%%%%%%%%%%%%%%%%%%%%%%%%%%%%%%%%%%%%%%%%%%%%%%%%%%%%%%%%%%%%%%%%%%%%%%%%%%%%%%%%%%%%%%%%%%%%%%%%%%%%%%%%%%%%%%%%%%%%%%%%%%%%%%%%%%%%%%%%%%%%%%%%%%%%%%%%%%%%%%%%%%%%%%%%%%%%%%%%%%%%%%%%%%%%%%%%%%%%%%%%%%%%%%%%%%%%%%%%%%%%%%%%%%%%%%%%%%%%%%%%%%%%%%%%%%%%%%%%%%%%%%%%%%%%%%%%%%%%%%%%%%%%%%%%%%%%%%%%%%%%%%%%%%%%%%%%%%%%%%%%%%%%%%%%%%%%%%%%%%%%%%%%%%%%%%%%%%%%%%%%%%%%%%%%%%%%%%%%%%%%%%%%%%%%%%%%%%%%%%%%%%%%%%%%%%%%%%%%%%%%%%%%%%%%%%%%%%%%%%%%%%%%%%%%%%%%%%%%

\section{The proof}
In this section we prove Theorem \ref{theorem} and Corollary \ref{corI}. At first, we recall in Section \ref{sectionbifurcation} the 
definitions of crossing forms and their application in bifurcation theory. Afterwards we compute the crossing forms for the equations \eqref{equationII} in Section \ref{computations}.

%%%%%%%%%%%%%%%%%%%%%%%%%%%%%%%%%%%%%%%%%%%%%%%%%%%%%%%%%%%%%%%%%%%%%%%%%%%%%%%%%%%%%%%%%%%%%%%%%%%%%%%%%%%%%%%%%%%%%%%%%%%%%%%%%%%%%%%%%%%%%%%%%%%%%%%%%%%%%%%%%%%%%%%%%%%%%%%%%%%%%%%%%%%%%%%%%%%%%%%%%%%%%%%%%%%%%%%%%%%%%%%%%%%%%%%%%%%%%%%%%%%%%%%%%%%%%%%%%%%%%%%%%%%%%%%%%%%%%%%%%%%%%%%%%%%%%%%%%%%%%%%%%%%%%%%%%%%%%%%%%%%%%%%%%%%%%%%%%%%%%%%%%%%%%%%%%%%%%%%%%%

\subsection{Crossing forms and bifurcation}\label{sectionbifurcation}
We follow in this section \cite{SFLPejsachowicz}. Let $H$ be an infinite dimensional separable real Hilbert space and let $I=[0,1]$ denote the unit 
interval. We consider $C^2$-functions $\psi:I\times H\rightarrow\mathbb{R}$ and assume throughout that $0\in H$ is a critical point of each 
functional $\psi_\lambda=\psi(\lambda,\cdot):H\rightarrow\mathbb{R}$. We call $\lambda_0$ a bifurcation point if any neighbourhood of $(\lambda_0,0)\in I\times H$ contains 
elements $(\lambda,u)$ such that $u\neq 0$ is a critical point of $\psi_\lambda$.\\ 
The Hessians $h_\lambda=D^2_0\psi_\lambda:H\rightarrow\mathbb{R}$ of $\psi_\lambda$ at $0\in H$ define a path of quadratic forms which we require henceforth to be 
continuously differentiable. Moreover, we assume that the Riesz representations of $h_\lambda$, $\lambda\in I$, with respect to the scalar product of $H$ are Fredholm operators.\\
We call $\lambda_0\in (0,1)$ a \textit{crossing} if $\ker h_{\lambda_0}\neq 0$. The \textit{crossing form} at $\lambda_0$ is defined by

\begin{align}\label{crossform}
\Gamma(h,\lambda_0):\ker h_{\lambda_0}\rightarrow\mathbb{R},\quad\Gamma(h,\lambda_0)[u]=\left(\frac{d}{d\lambda}\mid_{\lambda=\lambda_0}h\right)[u]
\end{align}
and a crossing $\lambda_0$ is said to be regular if $\Gamma(h,\lambda_0)$ is non-degenerate. It is easily seen that regular crossings are isolated and hence finite in number.\\
The proof of Theorem \ref{theorem} will be based on the theorem below which follows from Theorem 1 and Theorem 4.1 in \cite{SFLPejsachowicz} (cf. also \cite{Rabier}).

\begin{theorem}\label{bifurcation}
If $\lambda_0\in (0,1)$ is a regular crossing and 

\begin{align*}
\sgn\Gamma(h,\lambda_0)\neq 0,
\end{align*}
then $\lambda_0$ is a bifurcation point.
\end{theorem}

For the proof of Corollary \ref{corI} we need a further result from \cite{SFLPejsachowicz}. In what follows we assume that the quadratic forms $h_\lambda$, $\lambda\in I$, have finite Morse indices $M(h_\lambda)$. From Proposition 3.9 and Theorem 4.1 of \cite{SFLPejsachowicz} we obtain:

\begin{theorem}\label{Morse}
If all crossings of the path $h_\lambda$, $\lambda\in I$, are regular and $h_0,h_1$ are non-degenerate, then

\begin{align*}
M(h_0)-M(h_1)=\sum_{\lambda\in (0,1)}\sgn\Gamma(h,\lambda).
\end{align*}
\end{theorem}

%%%%%%%%%%%%%%%%%%%%%%%%%%%%%%%%%%%%%%%%%%%%%%%%%%%%%%%%%%%%%%%%%%%%%%%%%%%%%%%%%%%%%%%%%%%%%%%%%%%%%%%%%%%%%%%%%%%%%%%%%%%%%%%%%%%%%%%%%%%%%%%%%%%%%%%%%%%%%%%%%%%%%%%%%%%%%%%%%%%%%%%%%%%%%%%%%%%%%%%%%%%%%%%%%%%%%%%%%%%%%%%%%%%%%%%%%%%%%%%%%%%%%%%%%%%%%%%%%%%%%%%%%%%%%%%%%%%%%%%%%%%%%%%%%%%%%%%%%%%%%%%%%%%%%%%%%%%%%%%%%%%%%%%%%%%%%%%%%%%%%%%%%%%%%%%%%%%%%%%%%%
\subsection{Computation of the crossing forms}\label{computations}
We consider the family of boundary value problems \eqref{equation} for $0<r\leq 1$. After rescaling, equation \eqref{equation} on $\Omega_r$ is equivalent to 

\begin{equation}\label{equationres}
\left\{
\begin{aligned}
-\Delta u(x)+r^2\,g(r\cdot x,u(x))&=0,\quad x\in \Omega\\
u(x)&=0,\quad x\in\partial\Omega,
\end{aligned}
\right.
\end{equation}
and equation \eqref{equationII} reads as

\begin{equation}\label{equationIIres}
\left\{
\begin{aligned}
-\Delta u(x)+r^2\,f(r\cdot x)u(x)&=0,\quad x\in\Omega\\
u(x)&=0,\quad x\in\partial\Omega.
\end{aligned}
\right.
\end{equation}
Now consider the function $\psi:I\times H^1_0(\Omega)\rightarrow\mathbb{R}$ defined by 

\begin{align*}
\psi(r,u)=\frac{1}{2}\int_{\Omega}{\langle\nabla u,\nabla u\rangle\,dx}+r^2\int_{\Omega}{G(r\cdot x,u(x))\,dx},
\end{align*} 
where

\begin{align*}
G(x,t)=\int^t_0{g(x,\xi)\,d\xi}.
\end{align*}
Then $\psi$ is $C^2$ and the critical points of $\psi_r$, $r\in I$, are precisely the weak solutions of equation \eqref{equationres}. Moreover, $0\in H^1_0(\Omega)$ is a critical point for all $r$ and $r^\ast\in[0,1]$ is a bifurcation point of $\psi$ in the sense of Section \ref{sectionbifurcation} if and only if it is a bifurcation point for \eqref{equation}.\\ 
The Hessian of $\psi_r$ at the critical point $0$ is given by

\begin{align*}
h_r(u)=\int_{\Omega}{\langle\nabla u,\nabla u\rangle\,dx}+r^2\int_{\Omega}{f(r\cdot x)u^2\,dx},\quad u\in H^1_0(\Omega).
\end{align*}
Note that the kernel of $h_r$ consists of all classical solutions of the corresponding equation \eqref{equationIIres}. From the compactness of the embedding $H^1_0(\Omega)\hookrightarrow L^2(\Omega)$ it is easily seen that the Riesz representation of $h_r$ with respect to the usual scalar product of $H^1_0(\Omega)$ is a compact perturbation of the identity and hence a Fredholm operator. Moreover, the Morse index $M(h_r)$ is finite and so the function $\psi$ satisfies all assumptions of Section \ref{sectionbifurcation}.\\

\vspace*{0.2cm}

From the implicit function theorem applied to $\nabla\psi_r$, $r\in I$, we note at first that each bifurcation point of \eqref{equation} is a conjugate point of \eqref{equationII}. In order to show the remaining assertion of Theorem \ref{theorem} we now investigate the crossing form $\Gamma(h,r)$ at some crossing $r\in(0,1)$. We have by definition

\begin{align}\label{form}
\Gamma(h,r)[u]=\int_{\Omega}{\frac{d}{ds}\mid_{s=r}\left(s^2\,f(s\cdot x)\right)u^2(x)\,dx},\quad u\in\ker h_r.
\end{align}
Let now $0\neq u\in\ker h_r$ be given. Then $u$ is a classical solution of the equation \eqref{equationIIres}. We define for $s\in\mathbb{R}$ sufficiently small $u^s_r(x):=u(\frac{s}{r}x)$ and set

\begin{align}\label{diffu}
\dot u(x):=\frac{d}{ds}\mid_{s=r}u^s_r(x)=\frac{1}{r}\,\langle\nabla u(x),x\rangle.
\end{align}
Note that

\begin{align*}   
-\Delta u^s_r(x)+s^2\,f(s\cdot x)u^s_r(x)=0
\end{align*}
and, by differentiating this equation with respect to $s$ and evaluating at $s=r$, we get

\begin{align*}
-\Delta\dot u(x)+\frac{d}{ds}\mid_{s=r}\left(s^2\,f(s\cdot x)\right)u(x)+r^2\,f(r\cdot x)\dot u(x)=0.
\end{align*}
We multiply by $u$, integrate over $\Omega$ and conclude 

\begin{align*}
0&=-\int_{\Omega}{\Delta\dot u(x)u(x)\,dx}+\int_{\Omega}{\frac{d}{ds}\mid_{s=r}\left(s^2\,f(s\cdot x)\right)u(x)^2\,dx}\\
&+r^2\int_{\Omega}{f(r\cdot x)\dot u(x) u(x)\,dx}.
\end{align*}
Denoting by $\partial_nu(x)=\langle\nabla u(x),n(x)\rangle$, $x\in\partial\Omega$, the normal derivative to the boundary of $\Omega$, we obtain by applying Green's identity twice

\begin{align*}
0&=-\int_{\Omega}{\dot u(x)\Delta u(x)\,dx}+\int_{\partial\Omega}{(\partial_nu)\dot u\,dS}-\int_{\partial\Omega}{(\partial_n\dot u) u\,dS}\\
&+\int_{\Omega}{\frac{d}{ds}\mid_{s=r}\left(s^2\,f(s\cdot x)\right)u(x)^2\,dx}+r^2\int_{\Omega}{f(r\cdot x)\dot u(x) u(x)\,dx}.
\end{align*}
Now, since $u$ solves \eqref{equationIIres} and vanishes on $\partial\Omega$, we deduce by \eqref{form} and \eqref{diffu} that

\begin{align*}
\Gamma(h,r)[u]=-\int_{\partial \Omega}{(\partial_nu)\dot u\,dS}=-\frac{1}{r}\int_{\partial \Omega}{\langle\nabla u(x),n(x)\rangle\,\langle\nabla u(x),x\rangle\,dS}.
\end{align*}
Denoting by $x^T$ the tangential component of the vector $x\in\partial\Omega$, we have 

\begin{align*}
\langle\nabla u(x),x\rangle=\langle\nabla u(x),n(x)\rangle\langle x,n(x)\rangle+\langle\nabla u(x),x^T\rangle
\end{align*}
and hence

\begin{align*}
\Gamma(h,r)[u]&=-\frac{1}{r}\int_{\partial \Omega}{\langle\nabla u(x),n(x)\rangle^2\,\langle x,n(x)\rangle\,dS}-\frac{1}{r}\int_{\partial \Omega}{\langle\nabla u(x),n(x)\rangle\,\langle\nabla u(x),x^T\rangle\,dS}.
\end{align*}
It is easily seen that 

\begin{align*}
\langle\nabla u(x),n(x)\rangle\,\langle\nabla u(x),x^T\rangle=\diverg(u(x)(\partial_nu(x))x^T),\quad x\in\partial\Omega,
\end{align*}
and now we finally obtain by using Stokes' theorem on $\partial\Omega$

\begin{align*}
\Gamma(h,r)[u]&=-\frac{1}{r}\int_{\partial \Omega}{(\partial_nu(x))^2\langle x,n(x)\rangle\,dS}\leq 0,
\end{align*}
where we use that $\langle x,n(x)\rangle>0$ for all $x\in\partial\Omega$.\\
Moreover, if $\Gamma(h,r)[u]=0$, then $\partial_n u(x)=0$, $x\in\partial\Omega$, and we conclude that $u\equiv 0$ on $\Omega$ since $u$ solves the boundary value problem \eqref{equationIIres}.\\ 
Hence we have shown that $\Gamma(h,r)$ is negative definite. In particular, the crossing $r$ of $h$ is regular and 

\begin{align*}
\sgn\Gamma(h,r)=-\dim\ker h_r=-m(r).
\end{align*}
Now Theorem \ref{theorem} follows from Theorem \ref{bifurcation}, and Corollary \ref{corI} follows from Theorem \ref{Morse}.

\subsubsection*{Acknowledgement}
We want to thank Graham Cox for several helpful comments concerning the final step of our proof.

\vspace{1cm}
Alessandro Portaluri\\
Department of Agriculture, Forest and Food Sciences\\
Universit\`a degli studi di Torino\\
Via Leonardo da Vinci, 44\\
10095 Grugliasco (TO)\\
Italy\\
E-mail: alessandro.portaluri@unito.it

\vspace{1cm}
Nils Waterstraat\\
Dipartimento di Scienze Matematiche\\
Politecnico di Torino\\
Corso Duca degli Abruzzi, 24\\
10129 Torino\\
Italy\\
E-mail: waterstraat@daad-alumni.de

\end{document}